\documentclass[11pt]{article}

\usepackage{epsf,latexsym,amssymb}
\usepackage{graphicx}
\usepackage{amsmath}  %will give errors with \be \ee !!!
\newcommand{\comment}[1]{}
\newcommand{\edo}{\end{document}}

\newcommand{\R}{{\mathbb R}}  %ams bold
\newcommand{\Z}{{\mathbb Z}}  %ams bold
\newcommand{\abs}[1]{\left\vert #1 \right\vert}
\newcommand{\absn}[1]{\vert{#1}\vert} %smaller
\newcommand{\norm}[1]{\left\Vert #1 \right\Vert}
\newcommand{\norms}[1]{\left\vert #1 \right\vert}  %undecided what to sue for norm

\newtheorem{theorem}{Theorem}
\newtheorem{itlemma}{Lemma}[section] %number by section (set in \em by default)
\newtheorem{itproposition}[itlemma]{Proposition}
\newtheorem{itcorollary}[itlemma]{Corollary}
\newtheorem{itremark}[itlemma]{Remark}
\newtheorem{itdefinition}[itlemma]{Definition}
\newtheorem{itexample}[itlemma]{Example}

\newenvironment{lemma}{\begin{itlemma}\rm}{\end{itlemma}} %no-italics
\newenvironment{remark}{\begin{itremark}\rm}{\end{itremark}} %no-italics
\newenvironment{corollary}{\begin{itcorollary}\rm}{\end{itcorollary}}
\newenvironment{proposition}{\begin{itproposition}\rm}{\end{itproposition}}
\newenvironment{definition}{\begin{itdefinition}\rm}{\end{itdefinition}}
\newenvironment{example}{\begin{itexample}\rm}{\end{itexample}}

\newcommand{\be}[1]{\begin{equation}\label{#1}}
\newcommand{\ee}{\end{equation}}
\newcommand{\bl}[1]{\begin{lemma}\label{#1}}
\newcommand{\ble}[1]{\begin{lemmaex}\label{#1}}
\newcommand{\br}[1]{\begin{remark}\label{#1}}
\newcommand{\bt}[1]{\begin{theorem}\label{#1}}
\newcommand{\bd}[1]{\begin{definition}\label{#1}}
\newcommand{\bp}[1]{\begin{proposition}\label{#1}}
\newcommand{\bc}[1]{\begin{corollary}\label{#1}}
\newcommand{\bex}[1]{\begin{example}\label{#1}}
\newcommand{\ec}{\mybox\end{corollary}}
\newcommand{\eex}{\mybox\end{example}}
\newcommand{\eem}{\mybox\end{example}}
\newcommand{\el}{\mybox\end{lemma}}
\newcommand{\er}{\mybox\end{remark}}
\newcommand{\et}{\qed\end{theorem}}
\newcommand{\ed}{\mybox\end{definition}}
\newcommand{\ep}{\mybox\end{proposition}}
\newcommand{\epr}{\end{proof}}
\newcommand{\bpr}{\begin{proof}}

\newcommand{\ecs}{\end{corollary}}
\newcommand{\eexs}{\end{example}}
\newcommand{\els}{\end{lemma}}
\newcommand{\ers}{\end{remark}}
\newcommand{\ets}{\end{theorem}}
\newcommand{\eds}{\end{definition}}
\newcommand{\eps}{\end{proposition}}
\newcommand{\halmos}{\rule{1ex}{1.4ex}}
\newcommand{\qed}{\hfill \halmos} %put \qed at right margin
\newcommand{\mybox}{\hfill $\Box$} %put \qed at right margin (white square)
\newcommand{\beq}{\begin{eqnarray}}
\newcommand{\eeq}{\end{eqnarray}}
\newcommand{\beqn}{\begin{eqnarray*}}
\newcommand{\eeqn}{\end{eqnarray*}}
\newcommand{\bi}{\begin{itemize}}
\newcommand{\ei}{\end{itemize}}
\newcommand{\ben}{\begin{enumerate}}
\newcommand{\een}{\end{enumerate}}

\newcommand{\bes}[1]{\begin{subequations}\label{#1}\begin{eqnarray}}
\newcommand{\ees}[1]{\end{eqnarray}\end{subequations}}

\newcommand{\st}{\, | \,}

\newcommand{\dist}{{\rm dist}}

\newenvironment{proof}{\noindent {\em Proof}.\ }{\hspace*{\fill}$\halmos$\medskip}

\newcommand{\X}{\mathbb{X}}
\newcommand{\dX}{\partial\X}
\newcommand{\x}{\bar{x}}
\newcommand{\Rn}{\R^{n}}
\newcommand{\Rp}{\R_{\geq 0}}

\newcommand{\at}{\alpha }  % I had $\alpha _3, but I think not needed

\newcommand{\xo}{x^{\scriptscriptstyle 0}}

\newcommand{\xt}{x(t,\xo,u)}

\newcommand{\normi}[1]{\norm{#1}_{\infty }}
\newcommand{\kk}{{\mathcal K}}
\newcommand{\kl}{{\mathcal K}{\mathcal L}}
\newcommand{\ki}{{\mathcal K_\infty}}

\newcommand{\ggt}{\gamma }

\newcommand{\aunbar}{\alpha _1}
\newcommand{\aupbar}{\alpha _2}

\newcommand{\gV}{\nabla V}

\newcommand{\J}{{\mathcal J}}
\newcommand{\trace}{\mbox{trace}\,}
\newcommand{\Ex}{\mathbb{E}}

\newcommand{\ebar}{\bar{\varepsilon }}

\newcommand{\gVt}[1]{\gV(#1)^T}
\newcommand{\lbar}{{\bar \lambda }}
\usepackage{amsopn}

\DeclareMathOperator*{\argmin}{arg\,min}
\newcommand{\dVtilde}{{\widetilde \Delta }V}  % use it for V(F(x,q)) instead of V(f(x,u))

\usepackage[yyyymmdd]{datetime}
\def\mydate{\leavevmode\hbox{\short{\year}.\twodigits\month.\twodigits\day}}
\def\twodigits#1{\ifnum#1<10 0\fi\the#1}
\makeatletter
  \newcommand*\short[1]{\expandafter\@gobbletwo\number\numexpr#1\relax}
\makeatother

\title{Remarks on input to state stability of perturbed gradient flows,\\
motivated by model-free feedback control learning}
\author{Eduardo D. Sontag\\
\scalebox{.6}{\mydate:\currenttime}
}

\begin{document}

\maketitle

\begin{abstract}
  \noindent
Recent work on data-driven control and reinforcement learning has renewed
interest in a relative old field in control theory:~model-free optimal control
approaches which work directly with a cost function and do not rely upon
perfect knowledge of a system model. Instead, an ``oracle'' returns an
estimate of the cost associated to, for example, a proposed linear feedback
law to solve a linear-quadratic regulator problem.
%%%%%%%%%%%%%%%%%%%%%%%%%%%
This estimate, and an estimate of the gradient of the cost, might be obtained
by performing experiments on the physical system being controlled.  This
motivates in turn the analysis of steepest descent algorithms and their
associated gradient differential equations.
%%%%%%%%%%%%%%%%%%%%%%%%%%%
This note studies the effect of errors in the estimation of the gradient,
framed in the language of input to state stability, where the input represents
a perturbation from the true gradient.
%%%%%%%%%%%%%%%%%%%%%%%%%%%
Since one needs to study systems evolving on proper open subsets of Euclidean
space, a self-contained review of input to state stability definitions and
theorems for systems that evolve on such sets is included.
%%%%%%%%%%%%%%%%%%%%%%%%%%%
%Nothing surprising appears when generalizing from systems that evolve on
%Euclidean space, but it is hard to find a clean exposition in the literature.
%%%%%%%%%%%%%%%%%%%%%%%%%%%
The results are then applied to the study of noisy gradient systems, as well
as the associated steepest descent algorithms.
\end{abstract}

\section{Introduction}

Suppose that a function $V:\X\rightarrow \R$, defined on an open subset $\X$ of $\Rn$,
has a global minimum at a point $\x\in \X$, and that its gradient does not
vanish except at $x=\x$.  Under appropriate technical conditions, the
solutions of the gradient flow $\dot x=-\eta \gV(x)^T$ (where $\eta >0$ is a ``learning
rate'') will globally, and even exponentially, converge to $\x$ as $t\rightarrow \infty $.

In many data-driven applications, the gradient can be well-estimated numerically.
The combination of direct gradient estimation and gradient descent has
generated strong recent interest in control theory, and
specifically in Reinforcement Learning (RL) model-free control.  In order to
theoretically better understand the problem, several authors have studied an
archetypical control problem, the infinite-horizon Linear Quadratic Regulator
(LQR) problem.  Since the pioneering work of Kalman in the early 1960s, it has
been known that the solution of the LQR problem can be obtained explicitly
via a Riccati equation, and many computational packages do so very
efficiently.  Nonetheless, if the system being controlled is imperfectly
known, the function to be optimized is not known except through ``queries''
involving sampling and experimentation, and in that context direct methods
might be of interest. 
In any event, however, working on a well-understood problem like LQR serves to
understand properties of model-free approaches.

It turns out that when the LQR problem is formulated as an optimization over
a set of stabilizing feedback matrices, the loss function, while not convex,
satisfies strong convergence guarantees~\cite{2021mohammadi}.  The trick is to
employ a reparametrization for the LQR problem that allows solving an
associated strongly convex problem.  We refer the reader to~\cite{2021mohammadi} for
details.  Note that, in the LQR problem as just described, the open set $\X$
is a set of matrices.  Restricting the optimization dynamics to this open set
is essential for the approach to work.

In this note, we study a perturbed gradient system (superscript $T$ indicates transpose):
\be{eq:gradient}
\dot x(t) \;=\; -\eta \, \gV(x(t))^T\, + \,B(x(t))u(t) \,.
\ee
The additive term represents disturbances.
For example, if $B(x)$ is the constant matrix with rows $(1,0,\ldots ,0)$,
$(0,1,\ldots ,0)$, \ldots , $(0,0,\ldots ,1)$ then we have independent disturbances
$u_i$ acting on each coordinate.  
Without the additive term, this is a standard gradient descent flow.
For generality, we allow state-dependent perturbations (non-constant $B$).

The ``disturbance'' inputs might represent errors when numerically approximating
the gradient from data through two-point estimates as in~\cite{2021mohammadi},
or due to measurement noise.
The paper~\cite{gruene_iss_numerical} interprets the discretization error when
solving ODE's as a perturbation, and relates
asymptotic stability for dynamical systems to families of approximations,
specifically applying this to numerical one step schemes for ordinary
differential equations.

To quantify the effect of disturbances, we will use the
notion of input-to-state stability (ISS), introduced in~\cite{ISS89}
(see expositions in~\cite{04cime,mct,khalil,isidori}).
We will prove (under technical assumptions on $V$, mainly that $V$ blows up at
the boundary of $\X$, so that trajectories cannot escape the constraint set;
assumptions which hold in the motivating example from~\cite{2021mohammadi})
that the disturbed gradient system is ISS.
This implies that if the disturbances or errors  are bounded, small,
``eventually'' small, or convergent, the solutions of the system will inherit
the same properties, with well-controlled transient behavior.

The natural setting is that of differential equations that evolve in a nontrivial open
subset $\X$ of $\R^n$.  An example is a gradient system that uses a loss
function associated to a feedback matrix $K$ that is required to stabilize a
given linear system $\dot x=Ax+Bu$, in the sense that $A-BK$ is a Hurwitz matrix
(i.e., it has all its eigenvalues with negative real parts).  We can view
matrices of size $p\times q$ as elements of $\R^n$, $n=pq$.  Since eigenvalues
depend continuously on matrix entries (a standard fact, proved for example in
the linear algebra appendix in~\cite{mct}), the set $\X\subset \R^n$ of
stabilizing matrices (for a fixed system defined by $A$ and $B$) is open. 

The precise statement of the ISS result requires introducing appropriate notions
of stability and ``size'' of elements in open subsets.  We consider such
notions here; that material that we discuss, on ISS on open subsets, should be
of independent interest beyond the study of gradient systems, so we provide
detailed proofs of several facts about them for more general systems with
inputs. 

None of the results are surprising; nonetheless, it seems worth codifying the
basic theory with (almost) self-contained proofs.
Note that, if $\x\in \X$ is a globally asymptotically stable equilibrium, then
the set $\X$ must be diffeomorphic to $\Rn$ (this fact is actually used in the
proof of the converse ISS-Lyapunov theorem below).  So one could argue that,
up to this diffeomorphism, everything in this paper follows from the already
known results for systems in Euclidean space.  However, the diffeomorphism is
not \emph{a priori} known, and in any event, we wish to write everything in
the coordinates natural for the problem being studied.

In addition to studying the gradient system, we study the performance of
steepest descent, the discrete process in which a line search is performed,
iteratively minimizing a cost function in the direction of the gradient.
%%%%%%%%%%%%%%%%%%%%%%
Given a continuously function $V:\X\rightarrow \R$ to be minimized on an open subset
$\X\subseteq \Rn$, the steepest descent algorithm consists of the following procedure:
given any initial state $\xo$, one performs a line search in the negative
gradient direction so as to minimize $V(\xo-\lambda \gVt{\xo})$ over $\lambda \geq 0$; the
minimal point then defines a new point $x^1$, and one then iterates.  Observe
that this search only makes sense on a maximal interval such that the line
segment $\{\xo-\mu \gVt{\xo}, \mu \in [0,\lambda ]\}$ is included in $\X$ (so that one may
evaluate $V$ for increasing $\lambda $).  when the gradient is imperfectly
evaluated, the picture is further complicated by the fact that one in fact
moves in a direction $\xo-\lambda [\gVt{\xo})+B(x)u]$, for some unknown additive
``noise'' input vector $u$ (we include $B(x)$ to allow a state-dependence of
the input).  This gives an iteration that we write as
%%%%%%%%%%%%%%%%%%%%%%%%%%%%%%%%%%%%%%%%%%%%%%%%%%%%%%%%%%%%%%
$x^+ = x - \lambda \left[\gVt{x}+B(x)u\right]$.

It is in principle possible that even for a very small step one cannot
diminish the cost at all, and moreover one might even exit the set $X$
altogether for an input of large magnitude.  A trivial example of this is
provided by $\X=(-1,1)$, $B(x)=1$, and $V(x) = x^2/2$.  The perturbed steepest descent
procedure attempts to move to $x - \lambda  (x+u)$.  If we take any $x>0$ and any
$u< - x$ then for any step size $\lambda  >0$ the cost increases, which means that
the steepest descent procedure will be ``stuck'' at $x$.  Moreover, for
large $\lambda $ the expression $x - \lambda  (x+u)$ gives a result outside $\X$.  Of
course, this can be fixed if the magnitude of the input $u$ is ``not too
large'' compared to the state $x$.  Indeed, we will show that, under
reasonable technical assumptions, the steepest descent procedure is input to
state stable as a discrete-time system with respect to disturbances.

\section{Size functions on open subsets}

We start by introducing a notion of ``size'' that is well-suited to
quantifying global convergence to a given equilibrium, and which in particular
acts as a barrier function preventing escape from $\X$.  

\bd{definition:size}
Let $\X$ be an open subset of $\Rn$ and let $\x\in \X$.
We will say that
\[
\omega : \X \rightarrow  \R
\]
is a \emph{size function for} $(\X,\x)$ if $\omega $ is:
\ben
\item
continuous,
\item
positive definite with respect to $\x$, that is,
$\omega (\x)=0$ and 
$\omega (x)>0$ for all $x\in \X$, $x\not= \x$,
and
\item
proper, that is, for every real number $r\geq 0$, the sublevel set
$S_r:=\{x \st \omega (x)\leq r\}$ is a compact subset of $\X$.
\een
\eds

\br{remark:compact}
Observe that, since $\X$ is an open set, asking that $S_r$ is compact in the
induced topology of $\X$ is equivalent to asking that $S_r$ is compact as a
subset of $\R^n$.
\er

\br{remark:size}
Let us denote by $\norms{x}$ the standard Euclidean norm in $\Rn$ (any
other norm could be used as well).
When $\X=\Rn$, a natural choice of size is $\omega (x) = \norms{x-\x}$.
The notion that we introduce here is based on the beautiful paper of
Kurzweil~\cite{kurzweil}, which studied Lyapunov stability theory on open
sets, and is a particular case of ``measures'' in the sense of Lakshmikantham and
coauthors (see e.g.~\cite{lakhsmikantan}).  In~\cite{teel_praly2000,2017noroozi}, the concept is
called a ``proper indicator function'' (but we prefer not to use that term,
since ``indicator function'' is typically used for the characteristic function of
a set).
We remark that one could equally well define a size with respect to any closed
subset ${\cal A}$ of $\X$, simply asking that $\omega (x)=0$ if and only if $x\in {\cal A}$,
which is useful when studying convergence of solutions of differential
equations to non-point attractors.
Another point worth mentioning is that the definition of size function and
many of the results can equally well be formulated on a general differentiable
manifold $\X$; in a Riemannian manifold, one can take $\norms{x}$ as the
geodesic distance to $\x$, and all the functions of the form
$\omega (x)=\alpha (\norms{x})$ are size functions, when $\alpha $ is a function of class $\ki$ (defined below).
In that sense, the setup in this note is closely related to the work
in~\cite{2015_angeli_efimov}, in which a variant of ISS for systems evolving
in manifolds was considered. In that paper, the authors gave a definition that
relaxes the stability requirement for the unforced system; in this note,
instead, we study a notion which reduces to the usual one for systems in $\Rn$.
\er

The following elementary exercise in real analysis provides an intuitive characterization of size functions.  We denote by $\dX$ the boundary of the set $\X$ (which is empty if and only if $\X=\R^n$).

\bl{lemma:size}
The following two statements are equivalent for any function $\omega :\X\rightarrow \R$:
\bi
\item[(a)]
$\omega $ is a size function for $(\X,\x)$
\item[(b)]
$\omega $ is continuous, positive definite with respect to $\x$, and for every sequence $\{x_k\in \X,k\geq 1\}$,
\be{eq:diverge}
\mbox{if either $x_k\rightarrow \dX$ or $\norms{x}\rightarrow \infty $, necessarily $\omega (x_k)\rightarrow \infty $.}
\ee
\ei
\els

\bpr
We must show that property (\ref{eq:diverge}) is equivalent to compactness of every sublevel set $S_r$.

Suppose that property (\ref{eq:diverge}) is true, and pick any $r\geq 0$. By
Remark~\ref{remark:compact}, we need to prove that $S_r$ is closed and
bounded as a subset of $\Rn$.  

We first prove that $S_r$ is closed.  Suppose that a sequence $\{x_k\}$ in
$S_r$ is such that $x_k \rightarrow  x\in \Rn$ as $k\rightarrow \infty $. We must show that $x\in S_r$. 
There are two cases to consider: $x\not\in \X$ and $x\in \X$.
In the first case, being  the limit of elements in $\X$, necessarily $x\in \dX$.
Thus $x_k\rightarrow \dX$ and, by the assumed property, $\omega (x_k)\rightarrow \infty $, contradicting the fact that the sequence $\{\omega (x_k)\}$ is bounded (by $r$),  So this case cannot hold.
Thus $x\in \X$, so that $x\in S_r$ because $S_r$ is closed in the relative topology of $\X$.  
(More explicitly: by continuity of $\omega $, we have that $\omega (x_k)\rightarrow \omega (x)$, and hence $\omega (x) \leq  r$, so $x\in  S_r$.).

Next we prove that $S_r$ is bounded.  
Suppose by way of contradiction that there is a sequence $\{x_k\}$ in $S_r$ is
such that $\norms{x_k} \rightarrow  \infty $ as $k\rightarrow \infty $.  Again using the assumed property,
it follows that $\omega (x_k) \rightarrow  \infty $, contradicting 
that all $\omega (x_k)\leq r$.  Thus $S_r$ is bounded.

Conversely, suppose that $S_r$ is compact for every $r\geq 0$ and consider a
sequence $\{x_k\in \X,k\geq 1\}$.  Suppose first that $\norms{x_k} \rightarrow  \infty $.  We need
to prove that, for every $r > 0$, there is an integer $K$ so that $k > K \Rightarrow 
w(x_k) > r$. 
Suppose that this is not true, i.e.,
there is some $r$ and a subsequence $k_j\rightarrow \infty $ so that $x_{k_j}\in S_r$ for all
$k_j$. 
Replacing $\{x_k\}$ by this subsequence, we can then assume that $x_k\in S_r$
for all $k$, and still $\norms{x_k} \rightarrow  \infty $. Since $S_r$ is compact, there is a
convergent subsequence with its limit $x\in S_r$.    
This contradicts that $\norms{x_k} \rightarrow  \infty $.
Similarly, suppose that $x_k \rightarrow \dX$.  By contradiction, assume again that
there is some $r$ and a subsequence $k_j\rightarrow \infty $ so that $x_{k_j}\in S_r$ for all
$k_j$. 
Replacing ${x_k}$ by the subsequence, we can then assume that $x_k\in S_r$ for
all $k$, and still $x_k \rightarrow \dX$. 
By compactness, we can assume, taking a subsequence, that $x_k\rightarrow x\in S_r\subseteq \X$ for some $x$.
However, since $x_k \rightarrow \dX$ (because we have subsequences of a sequence converging to the boundary), this implies that $x\in \dX$.  We have a contradiction, because $\X$ and $\dX$ are disjoint subsets of $\Rn$.
This completes the proof.
\epr

Given any open set $\X\subseteq \Rn$ and any $\x\in \X$, there are many possible size functions for $(\X,\x)$.  As we remarked earlier, $\norms{x- \x}$ works when $\X=\Rn$. In general, we may use, for example:
\[
\omega (x) \;=\; \max\left\{
\norms{x-\x}\,,\,
\frac{1}{\dist(x,\dX)} - \frac{a}{\dist(\x,\dX)}
\right\}
\]
for any $a\geq 1$.  The case $a=2$ of this formula was given in~\cite{kurzweil},
and with that choice one has that $\omega (x)=\norms{x-\x}$ for all $x$ near $\x$.

\subsection{Comparing size functions}

In the same manner that any two norms on a finite dimensional space are equivalent, there is a notion of equivalence of size functions.

We denote by $\Rp$ the set of nonnegative real numbers.

Recall that ${\cal K}$ is the set of functions $\alpha :\Rp \rightarrow  \Rp$ that are continuous,
strictly increasing, and satisfy $\alpha (0)=0$, and ${\cal K}_\infty \subset{\cal K}$ is the subset
of unbounded functions, that is, $\alpha (r)\rightarrow \infty $ as $r\rightarrow \infty $.  The set ${\cal K}$ is
closed under sums, products, and compositions, as is the set ${\cal K}_\infty $.  Moreover,
functions in ${\cal K}$ are invertible, and $\alpha ^{-1}\in {\cal K}_\infty $, so ${\cal K}_\infty $ is a group
under composition (with identity element the map $\alpha (r)=r$). If $\alpha \in {\cal K}$, one
also says that ``$\alpha $ is of class ${\cal K}$'' and similarly for ${\cal K}_\infty $.  These
classes of functions have played a central role in dynamical systems since at
least the textbook by Hahn~\cite{hahn}, and were key in the development of
input to state stability notions in~\cite{ISS89}.  They have many other
useful properties, for example the weak subadditivity property
$\alpha (r+s)\leq \alpha (2r) + \alpha (2s)$; see for instance~\cite{mct,04cime}.  They allow
us to relate size functions. 
Observe that if $\omega $ is a size function and $\alpha \in \ki$, then $\alpha \circ \omega $ is also a
size function.

\bl{lemma:size_open_function}
Suppose that $\omega $ is a size function for $(\X,\x)$.  Then for each $\varepsilon >0$ there is
a $\delta >0$ such that
\[
   \omega (x) < \delta   \;\Rightarrow \; \norms{x-\x}< \varepsilon  \,.
\]
\els

\bpr
Since $\omega $ is continuous and $\omega (\x)=0$, there exists $\varepsilon _0>0$
such that $\abs{x-\x}<\varepsilon _0\Rightarrow \omega (x)<1$ (and the ball of radius $\varepsilon _0$ around
$\x$ is included in $\X$).
It follows that for any $x\in \X$ with $\abs{x-\x}\leq \varepsilon _0$, $\omega (x)\leq 1$.
Now pick any $\varepsilon >0$.  We let $\ebar:=\min\{\varepsilon ,\varepsilon _0\}$.
Consider the set
\[
C \,:=\; \{ x \st \norms{x-\x}\geq \ebar \mbox{ and } \omega (x)\leq 1 \}\,.
\]
The set is nonempty: pick any $x\in \X$ with $\abs{x-\x}=\varepsilon _0$; then
$\omega (x)\leq 1$ and also $\abs{x-\x}=\varepsilon _0\geq \ebar$.
The set $C$ is compact, because it is the intersection of a closed set and a compact set.
Also, $w(x)$ is nonzero in this set, because $w$ is positive definite.
Therefore there is a positive minimum of $w$ on the set $C$; we pick $\delta $ as this minimum, and thus
$x\in C \Rightarrow  \omega (x)\geq \delta $.
Without loss of generality, we will assume $\delta  < 1$ (otherwise, make $\delta $ smaller).
Now assume that $\omega (x)<\delta $.
This means that $x$ is not in $C$, so either $\omega (x)>1$ or $\norms{\x}<\ebar$.
However, $w(x)>1$ cannot happen,
because $\omega (x) < \delta  < 1$.
Therefore, $\norms{x-\x}<\ebar\leq \varepsilon $, as wanted.
\epr

\bp{proposition:compare}
Suppose that $\omega _1$ and $\omega _2$ are two size functions for $(\X,\x)$.  Then, there is some $\alpha \in {\cal K}_\infty $ such that 
\be{eq:compare}
\omega _1(x) \,\leq \, \alpha (w_2(x)) \mbox{ for all } x\in \X\,.
\ee
\eps

\bpr
Define
\[
\widetilde \alpha (r) \,:=\; \max_{\{x \st \omega _2(x)\leq r\}} \omega _1(x)\,.
\]
Since the set $\{x\st \omega _2(x)\leq r\}$ is compact, this maximum is well-defined.
Note that the inequality~(\ref{eq:compare}) holds.
Indeed, given any $x\in \X$, let $r:=\omega _2(x)$; then
$\omega _1(x)\leq \widetilde \alpha (r)=\widetilde \alpha (\omega _2(x))$, because $x$ belongs to the set over which we are maximizing.
Moreover, $\widetilde \alpha $ is nondecreasing (since as $r$ is larger, one takes a maximum over a larger set).
Also, $\widetilde \alpha (0)=0$ by positive definiteness of $\omega _1$ and $\omega _2$.  We prove next that $\widetilde \alpha $ is continuous at $0$.

Fix any $\varepsilon >0$.   We want to find a $\delta >0$ so that
   $r<\delta  \Rightarrow  \widetilde \alpha (r)<\varepsilon $.
From the definition of $\widetilde \alpha $, it is enough to find a $\delta $ such that, for each $r < \delta $:
\[
\omega _2(x) \leq  r \;\Rightarrow \; \omega _1(x) < \varepsilon /2\,.
\] 
Since $\omega _1$ is continuous and $\omega _1(\x)=0$, there is a $\delta _1>0$ such that
\[
\norms{x-\x} < \delta _1 \;\Rightarrow \; \omega _1(x) < \varepsilon /2 \,.
\]
By Lemma~\ref{lemma:size_open_function} applied to $\omega _2$ and $\varepsilon =\delta _1$, there is a $\delta >0$ such that
\[
   \omega _2(x) < \delta   \;\Rightarrow \; \norms{x-\x} < \delta _1\,.
\]
We conclude that:
\[
   \omega _2(x) < \delta  \;\Rightarrow \; \omega _1(x) < \varepsilon /2.
\]
Now assume $r<\delta $.
For any $x$ such that $\omega _2(x)\leq r$ , also $\omega _2(x) < \delta $, and hence $\omega _1(x) < \varepsilon $, as wanted.

So far we have a nondecreasing $\widetilde \alpha :\Rp\rightarrow \Rp$ which satisfies $\widetilde \alpha (0)=0$.
Such a function can be majorized by a class ${\cal K}_\infty $ function $\alpha $, i.e.\ $\widetilde \alpha (r)\leq \alpha (r)$ for all $r$,
which together with $\omega _1(x)\leq \widetilde \alpha (\omega _2(x))$ implies the estimate~(\ref{eq:compare}).
The construction of $\alpha $ is a standard exercise.
First majorize $\widetilde \alpha $ by a nondecreasing continuous function.  For example, pick a doubly infinite sequence of nonnegative numbers $r_k$, $k\in \Z$ so that $r_k\rightarrow 0$ as $k\rightarrow -\infty $ and $r_k\rightarrow \infty $ as $k\rightarrow +\infty $ and let $\alpha $ interpolate linearly the values $(r_k,\widetilde \alpha (r_{k+1}))$ (recall that $\widetilde \alpha $ is nondecreasing, so that the interpolation function is nondecreasing, and it clearly majorizes $\widetilde \alpha $).  This gives an $\alpha \in {\cal K}$.  Finally, add any ${\cal K}_\infty $ function to obtain an $\alpha $ of class ${\cal K}_\infty $.
\epr

\bc{corollary:comparisons}
Suppose that $\omega _1$ is a size function for $(\X,\x)$.  Let $\omega _2:\X\rightarrow \Rn$ be a
continuous function.  Then the following properties are
equivalent:
\bi
\item[(a)]
  $\omega _2$ is a size function for $(\X,\x)$;
\item[(b)]
  there exist functions $\alpha _1,\alpha _2\in \ki$ such that
\be{eq:compare-symmetric}
\alpha _1(\omega _1(x)) \,\leq \, \omega _2(x) \,\leq \, \alpha _2(w_1(x)) \;\;\mbox{ for all } x\in \X\,.
\ee
\ei
\ecs

\bpr
Suppose that $\omega _2$ is a size function for $(\X,\x)$.
By Proposition~\ref{proposition:compare}, there is an $\alpha \in \ki$ such that
$\omega _1(x)\leq \alpha (w_2(x))$ for all $x$.
Thus
$\alpha _1(\omega _1(x)) \leq  w_2(x)$
where $\alpha _1=\alpha ^{-1}$.
Applying again Proposition~\ref{proposition:compare}, but interchanging the
$\omega _i$'s, we have an $\alpha _2\in \ki$ such that
$\omega _2(x)\leq \alpha _2(w_1(x))$ for all $x$, so~(\ref{eq:compare-symmetric}) holds.

Conversely, suppose that~(\ref{eq:compare-symmetric}) holds.
Since $\omega _2(x)\geq \alpha _1(\omega _1(x))$ and $\alpha _1(\omega (x))>0$ for $x\not= \x$, it
follows that $\omega _2(x)>0$ for $x\not= \x$.  On the other hand,
$\omega _2(\x)\leq \alpha _2(w_1(\x))=0$, so $\omega _2$ is positive definite with respect to $\x$.
It remains to show that $\omega _2$ is proper.  Pick any $r\geq 0$ and consider
$S_r$.  This set is closed because $\omega _2$ is continuous.  On the other hand,
it is bounded because $S_r\subseteq \{x\st \omega _1(x) \leq  \alpha _1^{-1}(r)\}$ and the
latter set is compact because $\omega _1$ is proper.
\epr

\section{Systems with inputs}

From now on, assume given an open subset $\X\subseteq \R^n$, a point $\x\in \X$, and a
size function $\omega $ for $(\X,\x)$.
We consider here systems with $n$ state variables and $m$-dimensional inputs in the usual sense of control theory~\cite{mct}:
\[
\dot x(t) \; = \; f(x(t),u(t))
\]
(the argument ``$t$'' is often omitted, and dot indicates derivative with
respect to time).
The map
\[
f:\X\times \R^m\rightarrow \R^n
\]
is
assumed to be locally Lipschitz and $\x$ is an equilibrium when the input is
zero:
\[
f(\x,0)=0 \,.
\]
States $x(t)$ take values in $\X$, and inputs (also called ``controls'' or
``disturbances'' depending on the context)
are Lebesgue measurable essentially bounded maps
\[
u\,:\; [0,\infty )\rightarrow \R^m \,.
\]
We consider the sup norm of inputs:
\[
\normi{u}\,:=\; \mbox{ess}\sup_{t\geq 0}\norms{u(t)}
\]
where $\norms{u}$ is the Euclidean norm in $\R^m$ and ``ess sup'' denotes
essential supremum.  
%The norm of the restriction of an input to an interval $I$ is denoted by $\normi{u_I}$.

For each initial state $\xo$ and each input $u$, the solution of the initial
value problem with initial state $x(0)=\xo$ and input $u$ is denoted as
\[
\xt \;\in \;\X
\]
and is defined on some maximal interval
\[
  [0,t_{\max}(\xo,u))\,.
\]
%if $f$ is independent of inputs, simply as just as $\xtnou$.
%The \emph{zero-system} associated to $\dot x = f(x,u)$ is by definition
%the system with no inputs $\dot x=f(x,0)$.
%Euclidean norm is written as $\abs{x}$.

\br{remark:measurable}
For the sake of maximum generality, we allow inputs to be arbitrary (bounded)
measurable functions.  A technical issue is that measurable functions are in
reality equivalence classes of functions, equal only up to measure zero
subsets.  Solutions of the differential equation are absolutely continuous
functions and estimates over time have to be qualified by the phrase ``for
almost all $t$''.  We omit this qualification to make reading easier.  In any
event, for continuous inputs (which suffice for most applications) solutions are
continuously differentiable and there is no need for the qualifier.
\er

\subsection{Input to state stability}

The notion of input-to-state stability (ISS), introduced in~\cite{ISS89}
(see expositions in~\cite{04cime,mct,khalil,isidori})
provides a framework to describe stability features of the mapping
$(x(0),u(\cdot ))\mapsto x(\cdot )$
that sends initial states and input functions into solution trajectories.
Prominent among these features are that inputs that are bounded, small,
``eventually'' small, or convergent, should lead to states with the respective
property.  
In addition, ISS quantifies how initial states affect transient behavior.

The formal definition that we introduce, extended to open subsets, is as follows.
Recall that a function $\beta :[0,\infty )\times [0,\infty )\rightarrow [0,\infty )$ is said to be of class
$\kl$ if (1) for each fixed $t$, $\beta (s,t)$ as a function of $r$ is in class $\kk$ and (2) for each fixed $r$,
$\beta (r,t)$ decreases to zero as $t\rightarrow \infty $.

\bd{definition:ISS}
A system is \emph{input to state stable (ISS)} (on the open set $\X$ and with
respect to $\x$) if,
there exist functions $\beta \in \kl$ and $\ggt\in \ki$ so that the following property
holds:
for all inputs $u(\cdot )$ and all initial conditions $\xo\in \X$, the solution is
defined for all $t\geq 0$, that is, $t_{\max}(\xo,u)=+\infty $, and it satisfies the
estimate:
\[
\quad\quad\quad
\omega (\xt) \;\leq \;
\beta (\omega (\xo),t) \,+\, \ggt\left(\normi{u}\right)
\eqno(\mbox{ISS})
%\quad\quad\quad\quad\quad\quad(\mbox{ISS})
\]
for all $t\geq 0$.
\eds

Note that this definition is independent of the particular size function used
(although with $\beta $ and $\kappa $ functions that may change with $\omega $)
because of Proposition~\ref{proposition:compare}.
When $\X=\Rn$, since $\abs{x}$ is a size function, this becomes the usual
definition of ISS.

Since, in general, $\max\{a,b\}\leq a+b\leq \max\{2a,2b\}$,
one could restate the ISS condition in a slightly different manner, namely,
asking for the existence of some $\beta \in \kl$ and $\ggt\in \ki$ (in general
different from the ones in the ISS definition) such that
\[
\omega (\xt) \;\leq \;
\max\left\{\beta (\omega (\xo),t)\,,\, 
\ggt\left(\normi{u}\right)
\right\}
\]
holds for all solutions.

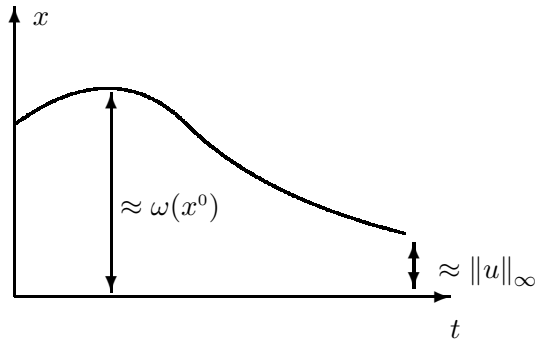
\begin{figure}[ht]%
%\begin{center}%
  \center
\setlength{\unitlength}{3000sp}%
%\begin{picture}(8424,2424)(2389,-3973)
\begin{picture}(8424,2424)(5000,-4100)
  \thicklines
\put(7201,-3961){\vector( 0, 1){2400}}
\put(7201,-3961){\vector( 1, 0){3600}}
\qbezier(7201,-2536)(8026,-1936)(8626,-2536)
\qbezier(8626,-2536)(9226,-3136)(10426,-3436)
\put(8000,-2275){\vector( 0, 1){  0}}
\put(8000,-2275){\vector( 0,-1){1650}}
%\put(8076,-3286){$\approx \vert \xo\vert$}
\put(8076,-3286){$\approx \omega (\xo)$}
\put(10801,-4300){$t$}
\put(7351,-1711){$x$}
\put(10501,-3511){\vector( 0, 1){  0}}
\put(10501,-3511){\vector( 0,-1){375}}
\put(10700,-3800){$\approx \normi{u}$}
\end{picture}
%\end{center}
\caption{ISS combines overshoot and asymptotic behavior}
\label{iss-fig1}
\end{figure}

Intuitively, the definition of ISS requires that, for $t$ large, the size of
the state must be bounded by some function of the sup norm, that is to say,
the maximum amplitude, of inputs, since $\beta (\omega (\xo),t)\rightarrow 0$ as $t\rightarrow \infty $.
On the other hand, the term $\beta (\omega (\xo),0)$
may dominate for small $t$, and this serves to quantify the magnitude of
the transient (overshoot) behavior as a function of the size of the initial
state $\xo$, see Figure~\ref{iss-fig1}.

For stable ($A$ having all eigenvalues with negative real part) linear systems
$\dot x=Ax+Bu$ evolving on $\X=\Rn$, the variation of parameters 
formula gives immediately the following inequality:
\[
\abs{x(t)} \;\leq  \; \beta (t)\abs{\xo} \,+\, \gamma  \normi{u}\,,
\]
where
\[
\beta (t) = \norm{e^{tA}} \; \rightarrow \;0 \quad \text{and}\quad
\gamma  = \norm{B} \int _0^{\infty }\norm{e^{sA}} ds \;<\;\infty 
\]
(here $\norm{\cdot }$ is induced operator norm).
This is a particular case of the ISS estimate, 
$\abs{x(t)}\leq \beta (\absn{\xo},t) +\gamma \left(\normi{u}\right)$,
with linear comparison functions.
Note that $\beta (t)\leq Ce^{-\lambda t}$ for some $C>0$ and some $\lambda >0$, so one has
exponential convergence when $u\equiv 0$.

\br{remark:sios}
We could think of a particular size function $\omega $ as an \emph{output} function
or ``observable'' $y=\omega (x)$ of the system $\dot x=f(x,u)$.  With this
interpretation, the definition is almost identical with that of
``state-independent input to output stability'' (SIIOS) given in~\cite{MR1754906}.
The paper~\cite{MR1754906} presents a large number of results relating SIIOS
to several other stability notions with respect to outputs.  However, the
interest in that paper is on \emph{non-proper} $\omega $, and $\X=\Rn$, since for
proper functions and $\X=\Rn$, SIIOS would simply coincide with ISS.
\er

\subsection{ISS-Lyapunov functions}

We now define ISS-Lyapunov functions on open sets.
We assume given a system $\dot x=f(x,u)$ as above.

\bd{definition:ISS-L}
A continuously differentiable $V:\X\rightarrow \R$ is said to be an
\emph{ISS-Lyapunov function} for $\dot x=f(x,u)$ if
\bi
\item[(a)]
$V-V(\x)$ is a size function for $(\X,\x)$, and
\item[(b)]
there exist functions $\at,\gamma \in \ki$ such that
\[
\dot  V(x,u) \; \leq \; -\at(\omega (x))\,+\,\gamma (\abs{u})
\quad\forall\, (x,u)\in \X\times \R^m
\eqno{(\mbox{L-ISS})}
\]
where
$\dot V:\X\times \R^m\rightarrow \R$ is the function:
\[
\dot V(x,u) := \nabla V(x).f(x,u)\,.
\]
\ei
\eds

The interpretation of $\dot V$ is given by the fact that, for any solution
$x(t)$ of $\dot x=f(x,u)$, the derivative $dV(x(t))/dt$ is $\dot V(x(t),u(t))$.

\br{remark:equiv-storage}
Property (a) in the definition of ISS-Lyapunov function is equivalent to
the existence of two functions $\alpha _i\in \ki$, $i=1,2$ such that
\be{eq:comparison-LISS}
\aunbar(\omega (x))\,\leq \,V(x)-V(\x)\,\leq \,\aupbar(\omega (x)) 
\quad\forall\,x\in \X\,.
\ee
This is an immediate application of Corollary~\ref{corollary:comparisons}.
Regarding property (b),
redefining $\alpha :=\at\circ \aupbar^{-1}\in \ki$, one also has an estimate in which,
instead of condition (\mbox{L-ISS}), one has the differential inequality:
\[
\dot  V(x,u) \; \leq \; -\alpha (V(x))-V(\x))\,+\,\gamma (\abs{u})
\quad\forall\, (x,u)\in \X\times \R^m \,.
\eqno{(\mbox{L-ISS'})}
\]
Conversely, suppose that (a) and (\mbox{L-ISS'}) hold.
Let $\alpha _1$ be as in~(\ref{eq:comparison-LISS}).
Then $\widetilde \alpha (\omega (x)) \leq  \alpha (V(x))-V(\x))$,
where $\widetilde \alpha :=\alpha \circ \aunbar$.
This $\widetilde \alpha $ gives an estimate of the form (\mbox{L-ISS}).
\er

\bt{theorem:iss-equiv}
A system is ISS if and only if it admits an %smooth
ISS-Lyapunov function.
\ets

The sufficiency of the ISS-Lyapunov condition is easy to show, and is entirely
analogous to the proof for $\X=\Rn$ in the original paper~\cite{ISS89}.
We sketch the details here, starting from an estimate L-ISS'.
Pick any solution $\xt$, and define
\[
v(t)\,:=\; V(\xt))-V(\x) \,.
\]
Note that $\dot v(t) = \dot V(x(t),u(t)) \leq  -\alpha (v(t))+\gamma (\abs{u(t)})$.
For any $t$, either $\alpha (v(t))\leq  2\gamma (\abs{u(t)})$ or
$\dot v(t)\leq  -\alpha (v(t))/2$.
From here, one deduces by a comparison theorem that
\[
v(t) \;\leq  \; \max\left\{\beta (v(0),t)\,,\,\at^{-1}(2\gamma (\normi{u}))\right\}
\quad\forall \, t \in [0,t_{\max}(\xo,u))\,,
\]
where the $\kl$ function
$\beta (s,t)$ is the solution $y(t)$ of the
initial value problem
\[
\dot  y=-\frac{1}{2}\at(y) %+\gamma (u)
\,,\quad y(0)=s \,.
\]
Using that
$v(0)=V(\xo)-V(\x)\leq \aupbar(\omega (\xo))$
and
$\omega (\xt) \leq  \alpha _1^{-1}(V(\xt)-V(\x))=\alpha _1^{-1}(v(t))$, we have
\beqn
\omega (\xt) &\leq &
\max\left\{\alpha _1^{-1}(\beta (\aupbar(\omega (\xo)),t))\,,\,\alpha _1^{-1}(\at^{-1}(2\gamma (\normi{u})))\right\}\\
&\leq &
\max\left\{\widetilde \beta (\omega (\xo),t)\,,\,\widetilde \gamma (\normi{u})\right\}
\quad\forall \, t \in [0,t_{\max}(\xo,u))\,,
\eeqn
with $~\beta \in \kl$ and $\widetilde \gamma \in \ki$.
It only remains to prove that $t_{\max}(\xo,u)=+\infty $.
To see this, note that, for any solution $\xt$, we have the bound
\[
\omega (\xt)\;\leq \; r := \max\left\{\widetilde \beta (\omega (\xo),0)\,,\,\widetilde \gamma (\normi{u})\right\}
\,.
\]
Therefore, $\xt\in S_r$ for all $t$ on the maximal interval of definition of the
solution. 
The set $S_r$ is compact (properness of size functions), so the solution is
defined for all $t\geq 0$ (see for example the ODE appendix in~\cite{mct}).

The converse part of the theorem follows by a reduction to the case $\X=\Rn$,
proved in~\cite{MR1372908,MR1325675}, which is basically a theorem about
Lyapunov functions for differential inclusions.
%The classical result of Massera~\cite{massera} for differential equations
%(with no inputs) becomes a special case.
Indeed, if a system is ISS, then the system with zero inputs $\dot x=f(x,0)$ has
$x=\x$ as an asymptotically stable point with domain of attraction all of $\X$
(this follows from the estimate $\omega (x)\leq \beta (\omega (0),t)$).
This implies that $\X$ is diffeomorphic to $\Rn$, see Theorem 2.2
in~\cite{wilson1967}, who obtains this as a simple corollary of the
Brown-Stallings Theorem.
(A proof of a simpler fact, that $\X$ must be contractible, is very easy; see
for example theorem 21 in~\cite{mct}.)
% got ref from Moulay's 2011 paper, who in turn also cites me :)
This means that under a diffeomorphism, we can apply the result for $\X=\Rn$,
and when transforming back, we obtain an L-ISS Lyapunov function.
\qed

\section{Application to gradient systems}

We assume given a pair $(\X,\x)$ and a size function $\omega $ for $(\X,\x)$.
We write the gradient of a function $V:\X\rightarrow \R$ as a row (co)vector
$\gV$, and its Euclidean norm as $\norms{\gV}$.
When $\gV$ is locally Lipschitz, the gradient flow has unique solutions
and if $\gV$ is globally Lipschitz, these solutions are automatically
defined for all $t\geq 0$.
(In the notations of Nesterov's book~\cite{2004nesterov}, the set of functions
$V$ for which $\gV$ has a uniform Lipschitz constant $L$ is denoted
$C_L^{1,1}(\X)$.  In our setup, solutions are defined for all
$t\geq 0$ even of $\gV$ is not assumed to be globally Lipschitz.)

\subsection{Proper loss functions}

\bd{definition:loss_function}
A continuously differentiable $V:\X\rightarrow \R$, with (locally) Lipschitz continuous gradient
$\gV$, will be said to be a \emph{proper loss function} if
\bi
\item[(a)]
$V-V(\x)$ is a size function for $(\X,\x)$, and
\item[(b)]
$\norms{\gV}$ is a size function for $(\X,\x)$.
\mybox
  \ei
\eds

\bl{lemma:loss_function_alternative}
Suppose given $V:\X\rightarrow \R$ continuously differentiable, with Lipschitz
continuous gradient $\gV$, such that $V-V(\x)$ is a size function for
$(\X,\x)$. Then these two properties are equivalent:
\bi
\item
  $V$ is a proper loss function,
\item  
  for some $\alpha \in \ki$,
\[
\alpha (\omega (x)) \;\leq \; \norms{\gV(x)}\qquad \mbox{for all } x\in \X\,.
\]
\ei
\els

\bpr
Suppose that $V$ is a proper loss function.
Since $\norms{\gV}$ is a size function for $(\X,\x)$, there is an $\alpha \in \ki$ as
claimed, by Corollary~\ref{corollary:comparisons} (take $\alpha =\alpha _1$).
Conversely, suppose that $\alpha (\omega (x)) \leq  \norms{\gV(x)}$ with $\alpha \in \ki$.
Then $\gV(\x)=0$, because $V$ has a (local and even global) minimum at $\x$.
For $x\not= \x$, $0 < \alpha (\omega (x)) \leq  \norms{\gV(x)}$, so $\norms{\gV(x)}$ is positive
definite.
For any $r\geq 0$, the set $\{x\st\norms{\gV(x)}\leq r\}$ is included in
$S_{\alpha ^{-1}(r)}$, so this set is bounded, and it is closed because $\gV$ is
continuous.  Thus $\gV$ is proper, and so it is a size function.
\epr

Applied with $\omega =V-V(\x)$, Lemma~\ref{lemma:loss_function_alternative}
together with the definition of size function says that
an equivalent way to define a proper loss function is to ask:
\bi
\item
  $V$ is continuously differentiable, with Lipschitz continuous gradient;
\item
  $V$ has a strict global minimum at $x=\x$;
\item
  $V$ is proper; and
\item
  there is some $\alpha \in \ki$ such that
\be{eq:compare-V-gradV}
\alpha (V(x)-V(\x)) \;\leq \; \norms{\gV(x)}^2 \qquad \mbox{for all } x\in \X\,.
\ee
\ei
We wrote $\norms{\gV(x)}^2$ instead of $\norms{\gV(x)}$ for convenience in what
follows; it makes no difference, since $(\alpha (\cdot ))^2$ is a $\ki$ function if and
only if $\alpha $ is.

In many problems one can directly obtain an
estimate as in~(\ref{eq:compare-V-gradV}), and this is useful for obtaining
explicit ISS stability rates.

\subsection{Gradient flow is ISS}

Fix a proper loss function $V$, a constant $\eta >0$ (the ``learning rate''),
and a locally Lipschitz and bounded mapping $B:\X\rightarrow \R^{n\times m}$.
We consider the %following
gradient system in~(\ref{eq:gradient}), repeated here for convenience:
\[
\dot x(t) \;=\; -\eta \, \gV(x(t))^T\, + \,B(x(t))u(t) \,.
\]

\bt{theorem:main}
If $V$ is a proper loss function, then system~(\ref{eq:gradient}) is ISS.
\ets

\bpr
We will prove that $V-V(\x)$ is an ISS-Lyapunov function for~(\ref{eq:gradient}).
Since $V-V(\x)$ is a size function, we need to show an estimate~(\mbox{L-ISS'}).
We have:
%for suitable $\alpha ,\gamma \in \ki$.
\beqn
\dot V(x,u) &=& - \eta \norms{\gV(x)}^2 \,+\, \gV(x)\,B(x)u\\
&=& - \eta \norms{\gV(x)}^2 + (\sqrt{\eta }\gV(x))(\sqrt{1/\eta } B(x)u)\\
&\leq & - \eta \norms{\gV(x)}^2 + \frac{\eta }{2}\norms{\gV(x)}^2 + \frac{1}{2\eta } \norms{B(x)u)}^2\\
&=& - \frac{\eta }{2}\norms{\gV(x)}^2 + \frac{1}{2\eta } \norms{B(x)u)}^2\\
&\leq & - \frac{\eta }{2}\alpha (V(x)-V(\x)) + \frac{C}{2\eta } \norms{u}^2\\
&=& -\widetilde \alpha ((V(x)-V(\x)) + \gamma (\norms{u})
\eeqn
where we used that, for row and column vectors in $\Rn$,
$\abs{vw} \leq  \norms{v}\norms{w} \leq  (1/2)(\norms{v}^2 + \norms{w}^2)$
(Cauchy-Schwarz inequality followed by $2ab\leq a^2+b^2$),
and the inequality $\alpha (V(x)-V(\x)) \leq  \norms{\gV(x)}^2$, and
where $C$ is an upper bound on $\norms{B}$,
and defined $\widetilde \alpha :=\frac{\eta }{2}\alpha \in \ki$ and $\gamma :=\frac{C}{2\eta }r^2\in \ki$.
So $V$ is an ISS-Lyapunov function, and thus the system~(\ref{eq:gradient}) is ISS.
\epr

In the particular case in which the estimate $\alpha (V(x)-V(\x))\leq \norms{\gV(x)}^2$
holds with a linear function $\alpha $, the proof of Theorem~\ref{theorem:iss-equiv}
provides a rate of decrease for $v(t)=V(x(t))-V(\x)$ which is exponential:
the function $\beta (r,t)$ has the form $e^{-\lambda t}r$ for some positive $\lambda $.

\subsection{An example: LQR problem}

The (infinite-horizon) LQR problem is one of the best-studied optimal control
problems.
Consider a time-invariant linear
system
\[
\dot x \;=\; Ax + B u
\]
and define the cost function:
\[
\J(\xo,u) \,:=\; \int_0^\infty  x^T (t) Q x(t) \, + \, u^T(t) R u (t) \,dt
\]
where $x(t)=\xt$.
Here $x(t)\in \Rn$, $u(t)\in \R^m$, 
$A\in \R^{n\times n}$ and $B\in \R^{n\times m}$ are matrices so that the pair $(A,B)$ is
controllable (or even just stabilizable or ``asymptotically controllable''),
which guarantees the finiteness of the objective function, and 
$Q\in \R^{n\times n}$, and $R\in \R^{m\times m}$ are positive definite.
The objective is to minimize $\J(\xo,u)$ over all measurable essentially bounded
control functions $u:[0,\infty )\rightarrow \R^m$, for any given $\xo$.

The unique optimal control is obtained by using the linear feedback law
$u(t)=-Kx(t)$, where $K = R^{-1}B^T\Pi $ and $\Pi $ is the unique positive
definite solution of the algebraic Riccati equation
\[
\Pi  B R^{-1}B^T\Pi  \,-\, A^T \Pi \,-\, \Pi  A \,-\, Q \;=\; 0
\]
(that is, $u(t) = -Kx(t)$), where $x$ solves $\dot x=(A-BK)x$ with $x(0)=\xo$),
and at this optimum value,
\[
\J(\xo,u) \;=\; (\xo)^T \Pi  \xo
\]
(see, for instance, Theorem 41 in~\cite{mct}).
The optimal feedback matrix $K= R^{-1}B^T\Pi $ stabilizes the
system, i.e., $A-BK$ is a Hurwitz matrix (all eigenvalues have negative part).

Since the optimal control is given by a linear feedback, one may pose the
simpler question of optimizing over all feedback matrices which belong to the
open set $\X:=\{K \st A-BK \mbox{ is Hurwitz}\}$.
In terms of $K$ and using $u=-Kx$,  one can introduce the loss function
\[
V_{\xo}(K)\,:=\; \int_0^\infty  x(t)^T Q x(t)  \; + \; (-Kx(t))^T \,R\, (-Kx(t)) \,dt
\]
where $x(t)$ solves $\dot x=(A-BK)x$, i.e., $x(t)=e^{(A-BK)t}\xo$, so that we can
also write
\[
V_{\xo}(K)\;=\; \int_0^\infty  x(t)^T (Q+K^TRK) x(t) \,dt
\;=\;
\trace\left( (Q+K^TRK) \int_0^\infty  x(t)\,x(t)^T\,dt \right)
\]
where we have used that for a scalar $a = \trace(a)$ and that
$\trace(UV)=\trace(VU)$.
To obtain a simpler problem, we assume now that the initial state is picked
distributed randomly according to some probability density in $\Rn$ (for
example, Gaussian) with covariance $\Sigma =\Ex[\xo (\xo)^T]$ and we wish to
minimize
%\be{loss-lqr}
\[
V(K)\,:=\; \Ex[V_{\xo}] \;=\; \trace((Q+K^TRK)P)
\]
%\ee
%where $N=Q+K^TRK$ and
with
\[
P\;=\; \Ex \left[\int_0^\infty  e^{tF}\xo (\xo)^T e^{tF^T}\,dt \right]
\;=\;
\int_0^\infty   e^{tF} \,\Sigma \, e^{tF^T}\,dt
\]
where $F=A-BK$.
It follows (see for instance Theorem 18 in~\cite{mct}) that
$P$ is the (unique) solution of the Lyapunov matrix equation
\be{eq:lyapunov-equation-linear}
(A-BK)P + P(A-BK)^T + \Sigma  \;=\; 0
\ee
In summary, one has to minimize the loss function
$V(K)=\trace((Q+K^TRK)P)$ where the positive definite matrix $P$
satisfies~(\ref{eq:lyapunov-equation-linear}) and $K\in \X$.
%This is a non-convex problem.
Since the solution $P$ of the linear system of
equations~(\ref{eq:lyapunov-equation-linear}) is a rational function of the
data (Cramer's rule), it follows that $V(K)$ is rational in the entries of the
matrix $K$, and hence $V$ is differentiable.
Although it is not generally convex, it has a unique global minimum at the
optimal $K= R^{-1}B^T\Pi $.  It is also known that it is a proper function,
see~\cite{1985toivonen}.
The gradient can be computed as follows (this is implicit in the computations
in~\cite{1970levine_athans,1985toivonen}, but see~\cite{1997rautert_sachs} for a clear
exposition):
\[
\gV(K)\;=\; 2 (RK - B^T L)P
\]
where $L$ is the unique positive definite matrix that satisfies
\[
(A-BK)^T L + L (A-BK) + Q + K^TRK \;=\; 0\,.
\]
For example, suppose that $n=m=1$, $a=q=r=\Sigma =1$. In this case
$\X=\{k\st bk>1\}$ and one obtains
\[
V(k) \;=\; \frac{k^2 + 1}{2(bk - 1)}
\]
and
\[
V'(k)\;=\; \frac{bk^2 - 2k - b}{2(bk - 1)^2}\,.
\]

In general, it can be shown, see~\cite{2021mohammadi}, that
$V$ is a proper loss function.
In fact, that reference shows that the 
Polyak-\L{}ojasiewicz condition~\cite{2016karimi}
\[
c_r (V(x)-V(\x)) \;\leq \; \norms{\gV(x)}^2
\]
holds on sublevel sets, for constants $c_r$, which implies that a lower
bounding $\alpha \in \ki$ exists.

\section{ISS and steepest descent}

From now on, we fix a pair $(\X,\x)$ and a size function $\omega $ for $(\X,\x)$.

We recall from the introductory discussion that we are interested in proving
that the steepest descent iteration $x^+ = x - \lambda [\gVt{x}+B(x)u]$, where $\lambda $
is picked at each step of the iteration so as to minimize the value $V(x^+)$,
is a discrete-time ISS system.

\subsection{Gradients of locally Lipschitz functions on $\X$}

Suppose given a continuously differentiable function $V:\X\rightarrow \R$ such that these
two properties hold:
\bi
\item[{\bf [SV]}]
$V-V(\x)$ is a size function for $(\X,\x)$,
\item[{\bf [LL]}]
$\gV$ is locally Lipschitz,
%\item[(c)]
%  $\gV(x)\not= 0$ for all $x\not= \x$.
\ei
We next review a couple of well-known facts about Lipschitz functions.

\br{remark:compact-lipschitz}
For each compact subset $K\subset X$, there some $L\geq 0$ such that the
one-sided Lipschitz estimate
\be{eq:1-L}
\left( \gV(y) - \gV(x)\right) (y-x) \;\leq \; L\abs{y-x}^2
\ee
holds for all $x,y\in K$.
Indeed, the function $\gV$ is Lipschitz on $K$, with some constant $L$ (start
locally and take finite subcovers),
% the argument does not need convexity!
so
\[
\abs{\left( \gV(y) - \gV(x)\right) (y-x)} \;\leq \;
\abs{\gV(y) - \gV(x)}\abs{y-x}  \;\leq \;
(L \abs{y-x})\abs{y-x} \;=\;
L\abs{y-x}^2
\]
by the Cauchy-Schwarz inequality and the Lipschitz property.
\er

\br{remark:nesterov}
Suppose that $x,y\in K$ and that $L$ is a one-sided Lipschitz constant as
in~(\ref{eq:1-L}) on the segment
\[
K \;=\; [x,y] \, :=\; \{z\st z=x+s(y-x),s\in [0,1]\}
\]
connecting $x$ and $y$.
Then
\be{eq:Lipschitz}
V(y)\;\leq \; V(x) + \gV(x)(y-x) + \frac{L}{2} \abs{y-x}^2 \,.
\ee
This is a standard fact, see e.g.~\cite{2004nesterov}.
The blanket assumption $\X=\Rn$ made there is not needed; since the proof is so
simple, we write it here.
Pick $x,y$, and consider the continuously differentiable function:
\[
g\;:\; [0,1]\rightarrow \R \;:\; s \mapsto  V(x+s(y-x)) \,.
\]
Then
\[
V(y)-V(x) - \gV(x)(y-x) \,=\,
g(1)-g(0) -g'(0) \,=\,
\int_0^1 g'(s)\,ds  -g'(0)\,=\,
\int_0^1 [g'(s)-g'(0)]\,ds
\]
where
\beqn
g'(s)-g'(0) &=&
\gV(x+s(y-x))(y-x) - \gV(x)(y-x)\\
&=&
\frac{1}{s} \left[ \gV(x+s(y-x)) - \gV(x)\right] (s(y-x)) \;\leq \;
sL\abs{y-x}^2
\eeqn
by~(\ref{eq:1-L}) when $s\not= 0$ (and this is trivial when $s=0$).
Therefore
\[
V(y)-V(x) - \gV(x)(y-x) \;\leq \; \int_0^1sL\abs{y-x}^2\,ds
\;=\; \frac{L}{2}\abs{y-x}^2\,,
\]
as desired.
\er

We wish to study the behavior of steepest descent when the gradient of $V$
is inaccurately estimated.

From now on we assume that $\gV$ is positive definite:
\bi
\item[{\bf [PD]}]
$\gV(x)\not= 0$ for all $x\not= \x$.
\ei
in addition to {\bf [SV]} and {\bf [LL]}.

\bl{lem:stay1}
Pick any $\xo\in \X$, $\xo\not= \x$, and let $L$ be a Lipschitz constant for $\gV$
on the compact set
\[
S\,:=\; \left\{x\in \X\st V(x)\leq V(\xo)\right\}
\]
(without loss of generality, $L>0$).
Pick any $q\in \Rn$ and write $p:=\gVt{\xo}\not= 0$.
Suppose that $\lambda >0$ has the property that
\[
%\xo-\mu (p+q)\in \X \;\mbox{and}\;
\xo-\mu (p+q)\in S \quad\mbox{for each}\; 0 \leq  \mu  \leq  \lambda \,.
\]
Then
\be{eq:descent}
V(\xo-\lambda (p+q)) - V(\xo) \;\leq \;
\left( -\lambda  + \frac{\lambda ^2L}{2}\right) \abs{p}^2 \,+\,
\frac{\lambda ^2L}{2}\abs{q}^2 \,+\,
\left( \lambda  + \lambda ^2L\right) \abs{p}\abs{q}\,.
\ee
\els

\bpr
Let $x=\xo$ and $y=\xo-\lambda (p+q)$.
The segment $[x,y]$ consists of points of the form
$\xo-\mu (p+q)$, with $0 \leq  \mu  \leq  \lambda $.
Therefore, we may apply the Lipschitz estimate~(\ref{eq:Lipschitz}), to obtain:
\[
V(\xo-\lambda (p+q)) - V(\xo) \;\leq \;
-\lambda p^T(p+q) \,+\,
\frac{\lambda ^2L}{2}\abs{p+q}^2 \,.
\]
Since
\[
\abs{p+q}^2 \;=\; \abs{p}^2+\abs{q}^2+2p^Tq \;\leq \; \abs{p}^2+\abs{q}^2+2\abs{p}\abs{q}
\]
and similarly $-p^Tq\leq \abs{p^Tq}\leq \abs{p}\abs{q}$, the
estimate~(\ref{eq:descent}) follows.
\epr

We have this immediate consequence:

\bc{cor:stay1}
Suppose that $q\leq c\abs{p}$ in Lemma~\ref{lem:stay1}.
Then,
\[
V(\xo-\lambda (p+q)) - V(\xo) \;\leq \;
\lambda \left[ (c-1) + \frac{\lambda L}{2} (c+1)^2 \right] \, \abs{p}^2\,.
\]
In particular, taking $c=1/2$ and $\lambda \leq \frac{2}{9L}$, then
$V(\xo-\lambda (p+q)) - V(\xo) \leq  -\frac{\lambda }{4}\abs{p}^2$.
\ec

\bl{lem:stay2}
Pick $\xo$, $L$, $q$, and $p$ as in Lemma~\ref{lem:stay1}.
Suppose that $\abs{q}\leq \frac{1}{2}\abs{p}$ and $\lambda =\frac{2}{9L}$.
Then $\xo-\mu (p+q)\in \X$ for each $0 \leq  \mu  \leq  \lambda $ and
\[
V(\xo-\lambda (p+q)) - V(\xo) \;\leq \; -\frac{1}{18 L}\abs{p}^2\,.
\]
\els

\bpr
Since $\X$ is an open set, $\xo-\mu (p+q)\in \X$ for all small $\mu >0$.
Also, since
\[
\left.\frac{d}{ds}\right|_{s=0} V(\xo-s(p+q))
\;=\; -p^T(p+q) 
\;=\; -\abs{p}^2 + p^Tq
\;\leq \; -\abs{p}^2 + \abs{p}\abs{q} 
\;\leq \; - \frac{1}{2}\abs{p}^2
\;<\; 0\,,
\]
there is some $\varepsilon >0$ such that
\[
V(\xo-\mu (p+q)) \;<\; V(\xo) \quad \mbox{for all}\; \mu \in (0,\varepsilon )\,.
\]
Suppose that there would exist some $\mu \in [0,\lambda ]$ such that
  $\xo-\mu (p+q)\not\in S$.
Since $S$ is compact and $\X$ is open, this would mean that there is some
$\mu \in [\varepsilon ,\lambda ]$ such that
$\xo-\mu (p+q)\in \X$ and $V(\xo-\mu (p+q))=V(\xo)$. %(*)
To apply Corollary~\ref{cor:stay1}, we need to see that this cannot happen.
Let
\[
\lambda _0\,:=\; \min \left\{\mu \in [\varepsilon ,\lambda ] \st V(\xo-\mu (p+q))=V(\xo)\right\} \;\geq \;\varepsilon   \;>\;0\,.
\]  
Since $V(\xo-\mu (p+q))\leq V(\xo)$ for all $\mu \in [0,\lambda _0]$, we may apply
Corollary~\ref{cor:stay1} to $\lambda _0$ to conclude that
\[
0 \;=\; V(\xo-\lambda (p+q)) - V(\xo) \;\leq \; -\frac{\lambda _0}{4}\abs{p}^2\,,
\]
which contradicts $\lambda _0>0$ and $p\not= 0$.
Thus the hypotheses of Lemma~\ref{lem:stay1} hold, and applying
Corollary~\ref{cor:stay1} to $\lambda $ we conclude that
\[
V(\xo-\lambda (p+q)) - V(\xo) \;\leq \; -\frac{\lambda }{4}\abs{p}^2 \;=\; -\frac{1}{18 L}\abs{p}^2
\]
as claimed.
\epr

\subsection{Line search in direction of steepest descent}

We continue with the assumptions {\bf [PD], [SV], [LL]} on $V$.
We next define a function
\[
F\,:\, \X\times \Rn \rightarrow  \X
\]
that will represent an individual steepest descent step when starting at a
point $\xo\in \X$ and the (transpose of the) gradient is estimated as $p+q$
where $p:=\gVt{\xo}$ and $q\in \Rn$ represents an additive noise.
%For any $\xo\in \X$, we let
%\[
%F(\xo,q) \,:=\; \xo
%\]
%when $p+q=0$, corresponding to an estimated zero gradient.
%We also let $F(\x,q)=\x$ if $q$ is arbitrary, as $V$ is already minimal there.
Now take any
%$\xo\not= \x$
$\xo\in \X$ and any $q\in \Rn$ such that $p+q\not= 0$.
Define
\[
\Lambda (\xo,q)\,:=\;
\left\{\lambda \geq 0 \st 
%\xo-\mu (p+q)\in \X \;\mbox{and}
V(\xo-\mu (p+q)) \leq  V(\xo) 
\;\mbox{for all}\;  \mu \in [0,\lambda ] \right\}\,.
\]
Note that $0\in \Lambda (\xo,q)$.
In the particular case $\xo=\x$, $V(\xo)$ is the unique minimizer of $V$,
so $\Lambda (\xo,q)=\{0\}$.  We consider $\xo\not= \x$ from now on.

The set $\Lambda (\xo,q)$ is compact.
It is bounded above: otherwise, it would be the case that $\xo-\lambda (p+q)\in \X$ and
$V(\xo-\lambda (p+q)) \leq  V(\xo)$ for all $\lambda \geq 0$; then since $V$ is proper, the set
of points $\xo-\lambda (p+q)$ is bounded, but this contradicts that
$\abs{\xo-\lambda (p+q)}\geq \abs{\abs{\xo}-\lambda \abs{p+q}}\rightarrow \infty $ as $\lambda \rightarrow \infty $ because
$\abs{p+q}\not= 0$.
It is also closed.
Indeed, suppose that $\lambda _k\rightarrow \lambda $, with $\lambda _k\in \Lambda (\xo,q)$.
Then $V(\xo-\lambda (p+q)) \leq  V(\xo)$, by continuity.
In addition, for each $\mu <\lambda $, there is some $k$ so that $\mu <\lambda _k$
so $V(\xo-\mu (p+q)) \leq  V(\xo)$, proving that $\lambda \in \Lambda (\xo,q)$.

Thus we may define
\[
\lbar_{\xo,q}\,:=\; \argmin_{\lambda \in \Lambda (\xo,q)} V(\xo-\lambda (p+q))
\]
where ``arg min'' means that we take the smallest $\lambda $ that achieves this
minimum value in the direction of $p+q$ when there is more than one.
We then define
\[
F(\xo,q) \,:=\; \xo-\lbar_{\xo,q}(p+q)\,.
\]
and $F(\xo,q):=\xo$ if $p+q=0$.
Note that $V(F(\xo,q))\leq V(\xo)$, because $0\in \Lambda (\xo,q)$ and we are minimizing.
In other words,
\[
\dVtilde(\xo,q)\,:=\; V(F(\xo,q)) - V(\xo) \;\leq \; 0 \quad \forall \, (\xo,q)
\]
and observe that $\dVtilde(x,q)=0$ if $V(\xo-\varepsilon (p+q))\geq V(\xo)$ for all small
$\varepsilon $. 
On the other hand, since
\[
\left.\frac{d}{ds}\right|_{s=0} V(\xo-sp) \;=\; -\abs{\gV(\xo)} < 0
\]
it follows that $\dVtilde(\xo,0)<0$ for all $\xo\not= \x$.

We next estimate $\dVtilde(\xo,q)$ for all $\abs{q}$ that are not ``too large''
    compared to $\abs{p}$.

Suppose that $L$ is any Lipschitz constant for $\gV$ on the set
$S=\left\{x\in \X\st V(x)\leq V(\xo)\right\}$, $\abs{q}\leq \frac{1}{2}\abs{p}$,
and $\lambda =\frac{2}{9L}$.
From Lemma~\ref{lem:stay2}, $\lambda \in \Lambda (\xo,q)$, so
$V(F(\xo,q))\leq V(\xo-\lambda (p+q))$ by definition of $\lbar_{\xo,q}$ as a minimizer.
Thus, again by the Lemma,
\be{eq:dt-LISSq}
\dVtilde(\xo,q)
\;=\; V(F(\xo,q)) - V(\xo)
\;\leq \; V(\xo-\lambda (p+q)) - V(\xo)
\;\leq \; -\frac{1}{18 L}\abs{\gV(\xo)}^2
\ee
where, recall, $L$ is a Lipschitz constant on $\{V(x)\leq V(\xo)\}$.

\subsection{Steepest descent with inputs}

We now consider a slightly more general setup as follows.
Let $B:\X\rightarrow \R^{n\times m}$ be a bounded mapping, and assume that the gradient error
at each iteration step is $q=B(x)u$, where $u\in \R^m$, so that
$\abs{q}\leq K\abs{u}$, where $K$ is an upper bound on the Euclidean induced norm
$\norm{B(x)}$, over all $x\in \X$.

We define the steepest descent algorithm, with inputs $u$, as the
discrete-time system defined by the following iteration function
$f:\X\times \R^m\rightarrow  \X$:
\[
x^+ \;=\; f(x,u) := F(x,B(x)u) \,.
\]
We define $\Delta V(x,u) := \dVtilde(x,B(x)u)$, that is
\[
\Delta V(x,u)\,:=\; V(f(x,u)) - V(x)\,.
\]
Since $\dVtilde(x,q)\leq 0$ for all $(x,q)$, also
$\Delta V(x,u)\leq 0$ for all $(x,u)$.
Obviously, we can also write
\[
\Delta V(x,u)\;=\; [V(f(x,u))-V(\x)] - [V(x)-V(\x)]\,.
\]
which exhibits $\Delta V$ as the change, in each steepest descent step, of the
``excess cost'' of $V$ compared to its minimum value $V(\x)$.

\subsection{Discrete-time ISS}

We now extend to open subsets the definition of input to state stability for
discrete time systems, which is completely analogous to that for continuous
time, see for
instance~\cite{99ifac,2001_automatica_dt_ISS_jiang_wang,2002_scl_jiang_wang}.
We consider discrete-time systems $x^+=f(x,u)$, where $f:\X\times \R^m\rightarrow \X$ is a
continuous function and $f(\x,0)=0$.

\bd{definition:DTISS}
The discrete-time system $x^+=f(x,u)$ is \emph{input to state stable (ISS)}
(on the open set $\X$ and with respect to $\x$) if
there exist functions $\beta \in \kl$ and $\ggt\in \ki$ so that the following property
holds:
for all input sequences $u=(u_0,u_1,\ldots )\in \ell^\infty _m$ and all initial conditions
$\xo\in \X$, the solution $\xo$ satisfies the estimate:
\[
\quad\quad\quad
\omega (\xt) \;\leq \;
\beta (\omega (\xo),t) \,+\, \ggt\left(\normi{u}\right)
\eqno(\mbox{ISS})
%\quad\quad\quad\quad\quad\quad(\mbox{ISS})
\]
for all $t=0,1,2,\ldots $.
\eds

Here $\normi{u}=\sum_{t=0}^\infty  \abs{u_t}$ and $\xt$ is obtained by solving
recursively $x_{t+1}=f(x_t,u_t)$.

There are several equivalent definitions of ISS-Lyapunov function for
discrete time systems.  We pick here the most convenient one for the current
application.

For any function $V:\X\rightarrow \R$, we denote $\Delta V(x,u):=V(f(x,u))-V(x)$.

\bd{definition:ISS-L-DT}
A continuous $V:\X\rightarrow \R$ is said to be an
\emph{ISS-Lyapunov function} for $x^+=f(x,u)$ if
\bi
\item[(a)]
$V-V(\x)$ is a size function for $(\X,\x)$, and
\item[(b)]
there are exist (i) a continuous and positive
definite function $\alpha $, and (ii) a function $\chi \in \ki$, such that:
\be{eq:dt-ISS}
\omega (x) \,\geq \, \chi (\abs{u})
\quad \Rightarrow  \quad
%\Delta V(x,u)\;\leq \; -\alpha (\omega (x))
\Delta V(x,u)\;\leq \; -\alpha (V(x)-V(\x))
\ee
for all $x\in \X$, $u\in \R^m$.
\ei
\eds

Equivalences among alternative ISS-Lyapunov function definitions, including a
condition of the type $\Delta V(x,u)\leq  -\at(\omega (x))+\gamma (\abs{u})$ for functions of class
$\ki$, are discussed
in Remark 3.3 of~\cite{2001_automatica_dt_ISS_jiang_wang}.
As with continuous-time systems, the existence of ISS-Lyapunov functions is
equivalent to the ISS property, see~\cite{99ifac,2002_scl_jiang_wang}. 
For completeness, and because of the interest in open subsets $\X$, we prove
the sufficiency below, appealing to some key technical lemmas
in~\cite{2002_scl_jiang_wang,2001_automatica_dt_ISS_jiang_wang}.

%The assumption that $\Delta V(x,u)\leq 0$ for all $(x,u)\in \X\times \Rn$ in the next result
%is made here in order to simplify the proof, but the theorem is still true
%without it, see \cite{2001_automatica_dt_ISS_jiang_wang}.

Let us write, for simplicity of notation, $W(x)=V(x)-V(\x)$.
As $\Delta V(x,u)=[V(f(x,u))-V(\x)] - [V(x)-V(\x)] = W(f(x,u))-W(x)$,
one can write~(\ref{eq:dt-ISS}) as:
%since $\Delta (x,u) = V(f(x,u))-V(x)$, one can rewrite
\[
\omega (x) \,\geq \, \chi (\abs{u})
\quad \Rightarrow  \quad
W(f(x,u)) \;\leq \; W(x) -\alpha (W(x))\,.
\]

\bt{theorem:DTiss-equiv}
If a discrete-time system admits an ISS-Lyapunov function $V$
%such that $\Delta V(x,u)\leq 0$ for all $(x,u)\in \X\times \Rn$,
then it is ISS.
\ets

\bpr
We first remark that one may redefine $V$, replacing it by a function of the
form $\rho (V(x))$ with $\rho \in \ki$, in such a manner that the
estimate~(\ref{eq:dt-ISS}) holds but now $\alpha \in \ki$ (and the redefined $V$ is
so that $V-V(\x)$ still a size function).
The argument is similar to the one given in~\cite{ISS89} for the continuous
time case, but it is more delicate, see the proof of Lemma 2.8
in~\cite{2002_scl_jiang_wang}.
Moreover, one may assume that $r\mapsto r-\alpha (r)$ is of class ${\cal K}$ (see Lemma B.1 in
\cite{2001_automatica_dt_ISS_jiang_wang}).
So, from now on, and redefining $V$ in this manner if needed,
we will assume that $\alpha $ satisfies these two properties.
Since $W$ is a size function, there is a $\pi \in \ki$ such that
$\omega (x)\geq \pi (W(x))$, and thus 
$W(x) \geq  \pi ^{-1}(\chi (\abs{u}))$ implies
$\omega (x)\geq \chi (\abs{u})$, so redefining $\chi $ as $\pi ^{-1}\circ \chi $ we can state the
ISS-Lyapunov property as:
\[
W(x) \,\geq \, \chi (\abs{u})
\quad \Rightarrow  \quad
%\Delta V(x,u)\;\leq \; -\alpha (\omega (x))
W(f(x,u)) \;\leq \; W(x) -\alpha (W(x))\,.
\]
Now let $\beta (r,t)$ be the solution of the scalar difference equation
\[
y_{t+1} \;=\; y_t-\alpha (y_t) \,,\quad y_0=r\geq 0\,.
\]
The property that $r\mapsto r-\alpha (r)$ is of class ${\cal K}$ implies $y_t\geq 0$ for all $t$,
and also that $y_t<y_t'$ implies $y_{t+1}<y_{t+1}'$ for any two solutions,
in other words, the iteration is monotone (it preserves order).
Thus the function $\beta $ is of class ${\cal K}$ on $r$.
Moreover, since $\alpha (y)\geq 0$, $y_{t+1}\leq y_t$,
%Moreover, $r-\alpha (r)$ of class ${\cal K}$ also implies that $r\geq \alpha (r)$, in other words,
iterates form a decreasing sequence.
Thus all solutions converge to zero as $t\rightarrow \infty $, since the only equilibrium $
y-\alpha (y)=y$ is at $y=0$.  So $\beta \in \kl$.

We introduce the following function $\gamma :[0,\infty )\rightarrow [0,\infty )$:
\[
\gamma (\mu )\,:=\; \max\left\{ W(f(x,u)) \st \abs{u}\leq \mu , W(x)\leq \chi (\mu )\right\}
\]
which is well-defined (the set over which we are maximizing is compact, and
$W(f(x,u))$ is continuous on $(x,u)$),
nondecreasing (the sets are larger as $\mu $ increases), and
satisfies $\gamma (0)=0$ (since $W(x)=0$ implies $x=\x$ and $f(\x,0)=\x$).
Note that this implication holds:
\[
W(x) \,\leq \, \chi (\abs{u})
\quad \Rightarrow  \quad
%\Delta V(x,u)\;\leq \; -\alpha (\omega (x))
W(f(x,u)) \;\leq \; \gamma (\abs{u}) \,.
\]
Replacing $\gamma $ by a larger function if needed, we may assume that $\gamma \in \ki$
and also that $\gamma (\mu )\geq \chi (\mu )$ for all $\mu $.
Consider the following sets:
\[
P_\mu  \,:=\; \left\{ x \st W(x)\leq \gamma (\mu ) \right\}\,.
\]
We claim that this set is forward invariant for inputs with $\norm{u}_\infty \leq \mu $.
Indeed, pick any $x\in P_\mu $ and any $u\in \R^m$ with $\norm{u}_\infty \leq \mu $.
If $W(x) \geq  \chi (\abs{u})$, then $W(f(x,u))\leq  W(x) \leq  \gamma (\mu )$, so $f(x,u)\in P_\mu $.
If instead $W(x) \leq  \chi (\abs{u})$, then
$W(f(x,u)) \leq  \gamma (\abs{u})\leq \gamma (\mu )$ as well.

Consider now any input $u$, any initial state $\xo$, and the corresponding
solution $\xt$ of $x^+=f(x,u)$.
Let $a_t:=W(\xt)$ for $t=0,1,\ldots $, and $\mu :=\norm{u}_\infty $.
We will compare this sequence to $y_t = \beta (W(\xo),t)=\beta (a_0,t)$.
Note that by definition $a_0=y_0$.

Consider first the case that $a_t\leq \gamma (\mu )$ for all $t$. Obviously in that case
$a_t\leq \max\{\beta (a_0,t),\gamma (\mu )\}$ for all $t$. 

Consider next the case that $a_t>\gamma (\mu )$ for some $t$.
Then either (i) $a_t>\gamma (\mu )$ for all $t$, or (ii) there is a $T\geq 0$ so that
$a_t\geq \gamma (\mu )$ for $t=0,\ldots ,T$ and $a_{T+1}\leq \gamma (\mu )$.
Suppose that $a_t\geq \gamma (\mu )$ for $t=0,\ldots ,T$.
From the ISS-Lyapunov property, we know that
\[
%a_t \,\geq \, \chi (M)
a_t \,\geq \, \gamma (\mu ) \quad \Rightarrow  \quad a_t \,\geq \, \chi (\mu ) \,\geq \, \chi (\abs{u_t})
\quad \Rightarrow  \quad
a_{t+1} \;\leq \; a_t -\alpha (a_t)\,.
\]
We claim that $a_t\leq y_t=\beta (a_0,t)$ for $t=0,\ldots ,T$.
This holds for $t=0$.
In general, if 
$\gamma (\mu )\leq a_t\leq y_t$ then $a_{t+1}\leq a_t -\alpha (a_t)\leq y_t-\alpha (y_t)=y_{t+1}$,
because $r-\alpha (r)$ is nondecreasing in $r$.
By induction, $a_t\leq y_t$ for $t=0,\ldots ,T$.
It cannot be that (i) holds, since $y_t\rightarrow 0$ as $t\rightarrow \infty $.
Thus (ii) holds.  Now the condition $a_{T+1}\leq \gamma (\mu )$ together with the
forward invariance of $P_\mu $ implies that $a_t\leq \gamma (\mu )$ for all $t>T$.

In summary, $a_t\leq \max\{\beta (a_0,t),\gamma (\mu )\}$ for all $t$, or
\[
W(\xt) \;\leq \; \max\{\beta (W(\xo),t),\gamma (\norm{u}_\infty )\}, \, t = 0,1,\ldots  \,.
\]
Let $\theta _i\in \ki$ be such that
$\omega (x)\leq \theta _1(W(x))$ and $W(x)\leq \theta _2(\omega (x))$.
Then
\[
\omega (\xt)
%comment out to avoid overfull
%\;\leq \; \theta _1(W(\xt))
\;\leq \; \max\{\theta _1(\beta (\theta _2(\omega (\xo)),t),\gamma (\norm{u}_\infty )\}
\;\leq \; \widetilde \beta (\omega (\xo),t) \,+\, \gamma (\norm{u}_\infty )
\]
with $\widetilde \beta (r,t) = \theta _1(\beta (\theta _2(r),t))$ is an ISS estimate.
\epr

\subsection{Application to steepest descent}

For each $r\geq 0$, we let
\[
L(r) \;=\; \mbox{a Lipschitz constant for} \, \gV \, \mbox{on the set}\,
\{x\st V(x)\leq r\}\,.
\]
Without loss of generality we may take $L$ as a continuous, nondecreasing, and
everywhere nonzero function.
Letting
\[
\theta (r)\,:=\; \frac{1}{18 L(r)}
\]
we conclude that:
\be{eq:dt-LISS1}
%\frac{1}{2K}
\abs{\gV(x)} \,\geq \, 2K \abs{u}
\quad \Rightarrow  \quad
\Delta V(x,u)\;\leq \; -\theta (V(x)) \abs{\gV(x)}^2 \,.
\ee
For $x=\x$, this is immediate since both sides vanish; for $x\not= \x$ it follows
from~(\ref{eq:dt-LISSq}).

From now on, we assume that:
\bi
\item[{\bf [SG]}]
  $\abs{\gV}$ is a size function for $(\X,\x)$.
\ei
This property implies {\bf [PD]}.

\bt{theorem:DTISS}
Suppose that {\bf [SV], [LL], [SG]} hold.
The system $x^+=f(x,u)$ is ISS, and $V$ is an ISS-Lyapunov function for it.
\ets

\bpr
We need to obtain an estimate as in~(\ref{eq:dt-ISS}).
Let $\theta $ be as in~(\ref{eq:dt-LISS1}).
Since $\abs{\gV}/(2K)$ is a size function, we may pick $\chi $ as any $\ki$
function with the property that $\chi (\abs{\gV(x)}/(2K))\geq \omega (x)$ for all $x$. 
Now, if the pair $(x,u)$ is such that $\omega (x)\geq  \chi (\abs{u})$, then
$\chi (\abs{\gV(x)}/(2K))\geq \omega (x)\geq \chi (\abs{u})$, and therefore
$\abs{\gV(x)}>2K\abs{u}$.   Thus we have the implication
\[
\omega (x) \,\geq \, \chi (\abs{u})
\quad \Rightarrow  \quad
\Delta V(x,u)\;\leq \; -\theta (V(x)) \abs{\gV(x)}^2 \,.
\]
Since both $\abs{\gV}^2$ and $V-V(\x)$ are size functions, there is some $\widetilde \alpha \in \ki$
such that $\abs{\gV(x)}^2\geq \widetilde \alpha (V(x)-V(\x))$,
from which
\[
-\theta (V(x))\abs{\gV(x)}^2 \;\leq \; -\theta (V(x))\, \widetilde \alpha (V(x)-V(\x))
\]
for all $x\in \X$.
Let
\[
\alpha (r)\; := \;\theta (r+V(\x)) \,\widetilde \alpha (r)\,.
\]
Since $\theta $ is continuous and everywhere positive, and $\widetilde \alpha \in \ki$, it follows
that $\alpha $ is continuous and positive definite.
We have
\[
\omega (x) \,\geq \, \chi (\abs{u})
\quad \Rightarrow  \quad
\Delta V(x,u)\;\leq \; - \alpha (V(x)-V(\x))
\]
and therefore $V$ is a discrete-time ISS Lyapunov function, as claimed.
\epr

We remark that in the special case that $\gV$ is globally Lipschitz, one can
take $L$, and hence also $\theta $, as a constant, so that $\alpha $ can be picked of
class $\ki$.

%\bp{prop:globally_L_steepest}
%Suppose that $\gV$ is globally Lipschitz.  Then, there are two positive
%constants $\chi $ and $\alpha $ such that
%\be{eq:globally_L_steepest}
%\abs{\gV(x)}  \,\geq \, \chi \abs{u}
%\quad \Rightarrow  \quad
%\Delta V(x,u)\;\leq \; - \alpha  \abs{\gV(x)}^2
%\ee
%for all $x\in \X$, $u\in \R^m$.
%\ep

\section{Discussion}

We have analyzed the ISS properties of continuous-time gradient descent on
open subsets of Euclidean space, as well as the ISS properties of the
associated discrete-time steepest descent algorithm.

The conditions that we impose, which generalize the Polyak-\L{}ojasiewicz
condition, have appeared in the recent literature in similar contexts.  For
example, in \cite{2021_arxiv_iss_gradient_suttner_dashkovskiy} one finds
extremum-seeking controllers based on gradient flows and an ISS property with
respect for disturbances, for an integrator and a kinematic
unicycle; $\X$ is a closed submanifold of $\Rn$.
The paper~\cite{2019_poveda_krstic_fixedtime_iss_extremum_seeking} studies the
gradient minimization of a function $V_q(x)$ on $\X=\Rn$, where a parameter $q$
represents time-varying uncertainty and an ISS property is established with
respect to the rate of change of $q$ (which is a notion called ``DISS'' in
\cite{DiffISS_IJRNC03}).  The
work~\cite{2020arxiv_bianchin_poveda_dallanese_gradient_iss_switched_systems}
solves an output regulation problem for switched linear dynamical systems,
with $\X=\Rn$, proving an ISS property for gradient flows with respect to
unknown disturbances acting on the plant.
In~\cite{ieee_tac_2018_cherukuri_et_al_convexity_saddle_point_dynamics}, the
authors also study a gradient flow and show ISS with respect to additive errors,
assuming strong convexity of the function to be minimized (in fact, a more
general ``convex-concave'' property), also with $\X=\Rn$.

It is worth remembering that ISS theory provides an overall conceptual view,
and is never the whole story.  To be useful in specific applications, good
estimates of the various gain functions are required.  An analogy is
Lyapunov-function analysis of nonlinear differential equations: while showing
stability is an important first step, in practice one wants quantifications of
overshoots, speed of convergence, and so on.  The brief discussion of the LQR
problem emphasizes that most of the actual work goes into establishing such
estimates, as is the case in the various works that we have cited.
Nonetheless, it seems useful to have a conceptual framework and ``roadmap''
that helps organize the overall abstract ideas.

Even at the conceptual level, there are many extensions still to be explored.
We mentioned the extension to Riemannian manifolds, which should be quite
straightforward.  More interestingly, finite- and fixed-time gradient
flows~\cite{2021garg,2020romero} may also be studied on open subsets and, most
importantly in the current context, from the point of view of finite and
fixed-time ISS in the sense of
e.g.~\cite{hong_jiang_feng_2008,2018lopezramirez,2019holloway_krstic}.

\section{Acknowledgments}

The author thanks Yuan Wang and Denis Efimov for a critical reading and
suggestions.  This research was supported in part by grants ONR N00014-21-1-2431
and AFOSR FA9550-21-1-0289.

\newpage
%\bibliographystyle{plain}
%\bibliography{2021_07_iss_gradient_systems}

\edo